\documentclass[10pt,twocolumn,twoside]{IEEEtran}
\usepackage{graphicx}
\usepackage{verbatim}
\usepackage{cite}
\usepackage[usenames]{color} 
\usepackage{soul,color}
\soulregister\cite7
\soulregister\ref7
\soulregister\pageref7
\soulregister\eqref7
\soulregister\split7

\usepackage{xcolor}

\usepackage{tabularx,booktabs}
\usepackage{makecell}
\newcolumntype{C}{>{\centering\arraybackslash}X}

\newtheorem{proposition}{Proposition}
\usepackage{amsfonts} 
\usepackage{lipsum}
\usepackage{multicol}
\usepackage{amsmath} 
\usepackage{bm}
\usepackage[utf8]{inputenc} 
\usepackage{authblk}
\usepackage{comment}
\ifCLASSOPTIONcompsoc
\usepackage[caption=false,font=normalsize,labelfon
t=sf,textfont=sf]{subfig}
\else
\usepackage[caption=false,font=footnotesize]{subfig}
\fi
\usepackage{empheq}
\usepackage{lipsum}
\usepackage{mathtools}
\usepackage{cuted}
\begin{document}
\title{A Two-stage Stochastic Model Considering Randomness of Demand in Variable Speed Limit and Boundary Flow Control}
\date{}
\author[a]{Hao Liu\thanks{hfl5376@psu.edu}}
\author[b]{Suyash Vishnoi\thanks{scvishnoi@utexas.edu}}
\author[b]{Christian Claudel\thanks{christian.claudel@utexas.edu}}
\affil[a]{Department of Civil and Environmental Engineering, The Pennsylvania State University, USA}
\affil[b]{Department of Civil, Architectural and Environmental Engineering, The University of Texas Austin, USA}
\maketitle
\begin{abstract}
The combination of boundary flow control and variable speed limit (VSL) is a widely used strategy for highway corridors to maintain safety and throughput when the capacity downstream is dropped. However, most proposed models assume fixed and known demand, which is not always true in reality. Following a recently proposed boundary flow and VSL control model that is based on the Lax-Hopf solution to the Lighthill-Whitham-Richards (LWR) model, this paper proposes a two-stage stochastic model to tackle the uncertainties in demand. Our model assumes the demand follows a discrete distribution and the speed limit can be chosen from a pre-defined set with reasonable values. The first stage tries to find the optimal inflow control, which is the here-and-now decision, by considering the speed limit control downstream, and the second stage is to determine the optimal speed limit after observing the real demand. Those controls are coordinated by using a rolling horizon scheme. The case study shows the proposed model outperforms three deterministic models, which use the minimum, mean value and maximum of the demand as model input, under various scenarios. With ensuring the maximum of throughput, the proposed model is demonstrated to be able to reduce the fluctuation of entry flow without a significant increase in blocked vehicles. 

\end{abstract}
\section{Introduction}
\IEEEPARstart{N}{owadays}, traffic congestion on highway networks is an inescapable and growing issue in cities across the world. It stems from the widespread increase in automobiles and results in a significant burden on both economy and environment. Due to the soaring cost and lack of available land, expanding the roadway network is no longer a sustainable way to mitigate congestion. Compared to building new roads, traffic flow control, which is an efficient alternative to combat the impacts of congestion, aims to manage the traffic in the best way with the available facility. There are assorted traffic flow control tools for highway networks, such as ramp metering \cite{papageorgiou2002freeway,gomes2006optimal,bellemans2006model}, variable speed limit \cite{carlson2011local, li2017reinforcement, hadiuzzaman2013cell}, toll lanes and the combination of these methods \cite{hegyi2005model, carlson2010optimal, vishnoi2021variable}.\\
Although traffic flow control has drawn considerable research efforts, it is still a challenging topic due to the uncertainties from a myriad of sources, which include both traffic flow model parameters such as road capacity\cite{hu2019multi, brilon2005reliability, geistefeldt2009comparative} and fundamental diagram \cite{wang2009speed, qu2017stochastic, li2012analysis}, and measured model inputs such as initial densities \cite{gazis1971line, liu2021robust, liu2018stochastic}, demand \cite{zhu2014accounting,li2017reinforcement, an2015robust} and turning ratios\cite{liu2021robustturning}. In most cases, such uncertainties come from either inevitable modeling errors or imperfect sensor measurements, and ignoring them can lead to poor control performance. For example, the ramp metering control highly depends on the downstream capacity, which could be a random variable. If it is overestimated in a deterministic control model, an excess of vehicles can be sent to the highway, and undesired consequence such as congestion can be induced. Otherwise, if the capacity is underestimated, a portion of the capacity can be wasted due to the insufficient inflow from the ramps. Therefore, it is essential for a control model to tackle the uncertainties mentioned above. The significance of handling such uncertainties has begun to attract researchers' interest. For instance, Schmitt and Lygeros \cite{schmitt2020convexity} proposed a robust control model for freeway networks, which minimizes the worst-case cost with uncertainties in both fundamental diagram and future demand.  Como et al. \cite{como2016convexity} proposed a traffic flow control method that is robust with respect to uncertainties in initial traffic volume and exogenous inflows. Liu et al. \cite{liu2021robustturning} proposed a robust control model to handle the uncertainties in turning ratios through distributionally robust chance constraints. \\
Although all the random parameters can impact control performance, it is impossible to consider all of them simultaneously in a single model. In fact, most robust traffic control models can only deal with one or two of them. The purpose of this paper is to propose a stochastic model to consider the uncertainties in the demand in a variable speed limit (VSL) and boundary control model. VSL is a widely studied topic since it can bring safety \cite{lee2006evaluation}, environment \cite{madireddy2011assessment} and mobility benefits \cite{carlson2010optimal}. Most VSL models are deterministic, and the uncertainties in demand, which play an important role on the control performance, is ignored. To the authors' best knowledge, studies on uncertainties in demand in VSL model are few and need more research efforts. Some existing examples include: Zhu and Ukkusuri \cite{zhu2014accounting} developed a link-based dynamic network loading (LDNL) model which computes the traffic dynamics between three parts of the original links- beginning, main and end. Then, they modeled the VSL problem as a Markov Decision Process (MDP) problem in which the demand and supply of the network are random and follow given distributions. The MDP is solved based on a reinforcement learning approach. Faisal et al. \cite{Alasiri2021RobustVS} proposed a robust VSL control model to consider the randomness in traffic dynamics by adding a noise term in the density evolution formula of a cell transmission model (CTM) \cite{daganzo1994cell, daganzo1995cell}. Then, they found a solution that guarantees that the traffic density asymptotically converges to the predetermined desired value. Li et al. \cite{li2017reinforcement} proposed a Q-learning (QL)-based VSL control strategy to relieve congestion at recurrent bottlenecks. The traffic dynamics are also simulated using CTM, and the model is validated in both stable and fluctuating demand scenarios.\\
In addition to the uncertainties, the importance of traffic flow modeling for control performance cannot be overstated either. As the only existing macroscopic traffic flow model, the Lighthill-Whitham-Richards (LWR) \cite{lighthill1955kinematic,richards1956shock} model has been used extensively since it was proposed. A common way to solve LWR is to approximate the derivatives with finite differences, e.g., CTM. More recently, Mazar{\'e} et al. \cite{mazare2011analytical} proposed a semi-analytical and grid free solution to the LWR model, the Lax-Hopf solution. This solution outperforms other solutions in two ways. First, unlike the models that discretize PDEs into ordinary differential equations (ODEs), e.g. \cite{carlson2011local, gugat2005optimal}, it does not require any discretization or approximation. Second, this solution can calculate the traffic state at any point directly from the initial and boundary conditions without any knowledge of prior event since it is grid free. Other analytical solutions like front tracking method \cite{henn2005wave} require full knowledge of prior events and may have exponential growth of waves in some situations as waves "bounce" back and forth from the boundary conditions, which results in low efficiency. Built upon the Lax-Hopf solution, different control models have been proposed: Li et al. \cite{li2014optimal} proposed a mixed integer linear program (MILP) to control the  boundary flow for highways; Liu et al. \cite{liu2021robust, liu2021robustturning} proposed robust control models to consider the uncertainties in initial densities and turning ratios, respectively; Vishnoi and Claudel \cite{vishnoi2021variable} proposed an MILP to control the boundary flows (including on-ramp flows and upstream mainlane flow) and VSL for highway networks.\\
However, \cite{vishnoi2021variable} assumes the travel demand is fixed and known. As mentioned before, demand is one of the main origins of uncertainties. As a follow-up research, this paper sets out to propose a two-stage stochastic model to address the uncertainties in demand. The proposed model assumes triangular FD for all links on the highway corridor. In this paper, the boundary and VSL control are implemented when the downstream capacity of the studied highway corridor is dropped due to certain reasons, such as crash and road work, to reduce the density while ensuring the maximization of throughput. The maximization of throughput only requires the demand on the corridor to be equal to the downstream capacity, which might block excessive vehicles at the entry node when the upstream demand is higher than the downstream capacity. To address this issue, we added a weighted penalty term for blocked vehicles at the entry node. In addition, the temporal variation in traffic flow has been demonstrated to affect safety on highways, see \cite{golob2004freeway, theofilatos2014review, yuan2021developing}. When traffic demand is random, a control based on a deterministic model may lead to significant inflow variation. Therefore, another penalty term for the flow fluctuation at the entry node is included as well. In the first stage, the feasible boundary control is defined. In the second stage, the real demand in a planning period is observed and the optimal VSL is solved. Then, a rolling horizon scheme is utilized to coordinate both controls.\\
The rest of this paper is organized as follows. Section \ref{sec:review} reviews the Lax-Hopf solution to the LWR model, and the constraints required for a VSL control model. Section \ref{sec:two-stage} develops the two-stage model to tackle the uncertainties in travel demand and the rolling horizon scheme to coordinate the boundary and VSL control. Section \ref{sec:casestudy} runs case studies and validates the model performance through comparing three deterministic models. We conclude in Section \ref{sec:summary}. \\
\section{Review of the VSL and Boundary Flow Control Model}\label{sec:review}
The VSL and boundary flow control model, which serves as the building block of this paper, is reviewed in this section. First, the Lax-Hopf solution according the triangular fundamental diagram (FD) and link-wise constraints, i.e. the compatibility conditions, are introduced. Next, the extra decision variables and constraints in the VSL and boudnary flow control model are presented. Then, we derive the two-stage model to address the uncertainties in demand in the next section.
\subsection{Lax-Hopf Solution}
This section briefly reviews the Lax-Hopf solution, which is a semi-analytical and grid-free solution, to the LWR model. The readers are referred to \cite{claudel2011convex, mazare2011analytical} regarding the details of the derivation.
\subsubsection{H-J PDE and Value Conditions}
The LWR model can be expressed as
\begin{equation}\label{eq:LWR}
\frac{\partial \rho(t,x)}{\partial t}+\frac{\partial \psi(\rho(t,x))}{\partial x}=0,
\end{equation}
where $x$ and $t$ are the spatial and temporal coordinate of a point, $\rho(t,x)$ is the corresponding density, $\psi$ denotes is defined as the fundamental diagram between flow and density. The proposed model is applicable for any concave and piecewise linear fundamental diagram. For non-linear FDs, the solution can be defined by replacing the original FDs with their piecewise linear approximations. This paper utilizes the triangular FD, which can be expressed as,
\begin{equation}\label{eq:fd}
\psi(\rho)= \max\{v_f\rho, w(\rho-\rho_m)\}\quad \rho \in [0, \rho_m]
\end{equation}
where $v_f$ is the free flow speed, $w$ is the congestion speed, $\rho_m$ is the jam density, where the flow is zero due to the total congestion. $\rho_c=-\frac{\rho_mw}{v_f-w}$ is the critical density where the flow reaches its capacity. By integrating the LWR PDE in space, another traffic flow model, H-J PDE, can be expressed as
\begin{equation}\label{eq:HJ}
\frac{\partial M(t,x)}{\partial t}-\psi(-\frac{\partial M(t,x)}{\partial x})=0,
\end{equation}
where $M(t,x)$, known as the Moskowitz function \cite{moskowitz1965discussion}, denotes the index of the vehicles at point $(t,x)$.\\
The spatial domain $[\xi, \chi]$ is divided evenly into $k_{max}$ segments and the time domain $[0, t_{max}]$ is divided evenly into $n_{max}$ segments. Let $K=\{1,...,k_{max}\}$ and $N=\{1,...,n_{max}\}$. Assuming the initial density in each spatial segment and the flow in each time step are constant, the piecewise affine initial condition $M_k(t, x)$, upstream boundary condition $\gamma_n(t, x)$, and downstream boundary condition $\beta_n(t, x)$ are defined as
\begin{equation}\label{eq:ic}
M_k(t, x)=
\begin{cases}
-\sum_{i=1}^{k-1}\rho(i)X\\
-\rho(k)(x-(k-1)X), & \text{if $t=0$}\\
& \text{and $x\in[(k-1)X, kX]$}\\

+\infty, & \text{otherwise}

\end{cases}
\end{equation}

\begin{equation}\label{eq:uc}
\gamma_n(t, x)=
\begin{cases}
\sum_{i=1}^{n-1}q_{\text{in}}(i)T\\
+q_{\text{in}}(n)(t-(n-1)T), & \text{if $x=\xi$}\\
& \text{and $t\in[(n-1)T, nT]$}\\

+\infty, & \text{otherwise}

\end{cases}
\end{equation}

\begin{equation}\label{dc}
\beta_n(t, x)=
\begin{cases}
\sum_{i=1}^{n-1}q_{\text{out}}(i)T\\
+q_{\text{out}}(n)(t-(n-1)T)\\
-\sum_{k=1}^{k_{max}}\rho(k)X, & \text{if $x=\chi$}\\
& \text{and $t\in[(n-1)T, nT]$}\\

+\infty, & \text{otherwise}

\end{cases}
\end{equation}
where $X$ and $T$ are the spatial segment length and time step size, respectively, $\rho(i)$ is the initial density for the $i$th spatial segment and $q_{\text{in}}(i)$ and $q_{\text{out}}(i)$ are the inflow and outflow for the $i$th time step, respectively. For simplicity, the initial and boundary conditions are collectively referred to as value conditions (VCs).
\subsubsection{Semi-analytical and Grid-free solutions}
Given the VCs defined in the previous section, the Moskowitz solution at point $(t,x)$ from each VC can be computed by using the Lax-Hopf formula.\\
\begin{proposition}
	\textit{(Lax-Hopf Formula)}: Let $\psi(\cdot)$ be a concave and continuous Hamiltonian, and let $c(\cdot, \cdot)$ be a value condition. The B-J/F solution $M_c(\cdot, \cdot)$ to \eqref{eq:HJ} associated with $c(\cdot, \cdot)$ is defined \cite{aubin2008dirichlet, claudel2010lax, claudel2010lax2} by
	\begin{equation}\label{eq:lax-hopf}
	M_c(t,x)=\inf_{(u,T)\in(\varphi^*)\times R_+}(c(t-T, x+Tu)+T\varphi^*(u))
	\end{equation}
	where $\mathbf{\varphi^*}(\cdot)$ is the Legendre-Fenchel transform of an upper semicontinuous Hamiltonian $\psi(\cdot)$, which is given by,
	\begin{equation}\label{eq:lf}
	\varphi^*(u):= \sup_{p\in Dom(\psi)}[p\cdot u+\psi(p)]
	\end{equation}
\end{proposition}
By substituting the triangular FD \eqref{eq:fd} and the VCs \eqref{eq:ic}-\eqref{dc} into \eqref{eq:lax-hopf}-\eqref{eq:lf}, the corresponding Moskowitz solution can be expressed as
\scriptsize
\begin{subequations}\label{eq:mmk}
	\begin{empheq}[left={M_{M_k}(t, x)=\empheqlbrace\,}]{alignat=2}
	&+\infty, && \quad \text{if $x\le (k-1)X+tw$}\\
	& && \quad \text{or $x\ge kX+v_ft$}\nonumber\\
	&-\sum_{i=1}^{k-1}\rho(i)X+\rho(k) && \quad \text{if $x\ge (k-1)X+v_ft$}\\
	&(tv_f+(k-1)X-x), && \quad \text{and $x\le kX+v_ft$}\nonumber\\
	& && \quad \text{and $\rho(k)\le \rho_c$}\nonumber\\
	&-\sum_{i=1}^{k-1}\rho(i)X+\rho_c && \quad \text{if $x\le (k-1)X+v_ft$}\label{eq:example} \\
	&(tv_f+(k-1)X-x), && \quad \text{and $x\ge (k-1)X+tw$}\nonumber\\
	& && \quad \text{and $\rho(k)\le \rho_c$}\nonumber\\
	&-\sum_{i=1}^{k-1}\rho(i)X+\rho(k) && \quad \text{if $x\le kX+tw$}\\
	&(tw+(k-1)X-x)\nonumber && \quad \text{and $x\ge (k-1)X+tw$}\nonumber\\
	&-\rho_mtw, && \quad \text{and $\rho(k)\ge \rho_c$}\nonumber\\
	&-\sum_{i=1}^{k}\rho(i)X && \quad \text{if $x\le kX+tv_f$}\\
	&+\rho_c(tw+kX-x)\nonumber &&\quad \text{and $x\ge kX+tw$}\nonumber\\
	&-\rho_mtw, &&\quad \text{and $\rho(k)\ge \rho_c$}\nonumber
	\end{empheq}
\end{subequations}\\
\begin{subequations}\label{eq:mgamma}
	\begin{empheq}[left={M_{\gamma_n}(t, x)=\empheqlbrace\,}]{alignat=2}
	& +\infty, && \quad \text{if $t\le (n-1)T+\frac{x-\xi}{v_f}$}\\
	&\sum_{i=1}^{n-1}q_{\text{in}}(i)T+q_{\text{in}}(n) && \quad \text{if $t\ge (n-1)T+\frac{x-\xi}{v_f}$}\\
	&(t-\frac{x-\xi}{v_f}-(n-1)T), && \quad \text{and $t\le nT+\frac{x-\xi}{v_f}$}\nonumber\\
	&\sum_{i=1}^{n}q_{\text{in}}(i)T+\rho_cv_f && \quad \text{otherwise}\\
	&(t-\frac{x-\xi}{v_f}-nT), \nonumber
	\end{empheq}
\end{subequations}\\
\begin{subequations}\label{eq:mbeta}
	\begin{empheq}[left={M_{\beta_n}(t, x)=\empheqlbrace\,}]{alignat=2}
	&+\infty, && \quad \text{if $t\le (n-1)T+\frac{x-\chi}{w}$}\\
	&-\sum_{k=1}^{k_{max}}\rho(k)X+ && \quad \text{if $t\ge (n-1)T+\frac{x-\chi}{w}$}\\
	&\sum_{i=1}^{n-1}q_{\text{out}}(i)T+&&\quad \text{and $t\le nT+\frac{x-\chi}{w}$}\nonumber\\
	&q_{\text{out}}(n)(t-\frac{x-\chi}{w}, && \nonumber\\
	&-(n-1)T)-\nonumber\\
	&\rho_m(x-\chi),\nonumber\\
	&-\sum_{k=1}^{k_{max}}\rho(k)X && \quad \text{otherwise}\\
	&+\sum_{i=1}^{n}q_{\text{out}}(i)T+\nonumber\\
	&\rho_cv_f(t-nT-\frac{x-\chi}{v_f}), \nonumber
	\end{empheq}
\end{subequations}\\
\normalsize
\hrulefill
\vspace*{4pt}\\
\normalsize
\eqref{eq:mmk}-\eqref{eq:mbeta} shows that at a given point $(t,x)$, each VC generates a Moskowitz solution. According to \cite{claudel2011convex}, these solutions need to satisfy the inf-morphism property: the real Moskowitz solution is equal to the minimum among the ones from all VCs. Therefore, to compute the Moskowitz solution at any given point, we can use the Lax-Hopf formula to compute the Moskowitz solution from each VC directly and find the minimum without any knowledge about the prior events.\\
\subsubsection{Linkwise Constraints}
According to the inf-morphism property, to ensure the Moskowitz solutions at the domain of VCs are equal to the VCs, i.e., to make the VCs compatible with each other, we need the Moskowitz solutions at the domain of VCs are larger than or equal to the VCs. 
\begin{proposition}
	\textit{(Compatibility Conditions)}: Use the value condition \boldsymbol${c}(t, x)$ and the corresponding solution in \textit{Proposition 2}. The equality $\forall (t, x)\in Dom(\boldsymbol{c}), \boldsymbol{M_c}(t,x)=\boldsymbol{c}(t, x)$ is valid if and only if the inequalities below are satisfied,
	\begin{equation}
	\boldsymbol{M_{c_j}}(t, x)\ge \boldsymbol{c}_i(t, x), \quad \forall (t, x)\in Dom(\boldsymbol{c}_i), \forall(i, j)\in J^2
	\end{equation}
\end{proposition}
In detail, these constraints can be expanded as \cite{li2014optimal,canepa2012exact},
\begin{equation}\label{eq:ineq1}
\begin{cases}
M_{M_k}(0, x_p)\ge M_p(0, x_p) & \text{$\forall(k,p)\in K^2$}\\
M_{M_k}(pT, \chi)\ge \beta_p(pT, \chi) & \text{$\forall k\in K, 
	\quad \forall p\in N$}\\
M_{M_k}(\frac{\chi-x_{k}}{v_f}, \chi)\ge
\beta_p(\frac{\chi-x_{k}}{v_f}, \chi) & \text{$\forall k\in K, 
	\quad \forall p\in N$}\\
\hfill \text{s.t.}\quad \frac{\chi-x_{k}}{v_f}&\in[(p-1)T,pT]\\
M_{M_k}(pT, \xi)\ge \gamma_p(pT, \xi) & \text{$\forall k\in K, 
	\quad \forall p\in N$}\\
M_{M_k}(\frac{\xi-x_{k-1}}{w}, \xi)\ge 
\gamma_p(\frac{\xi-x_{k-1}}{w}, \xi) & \text{$\forall k\in K, 
	\quad \forall p\in N$}\\
\hfill \text{s.t.}\quad \frac{\xi-x_{k-1}}{w}&\in[(p-1)T,pT]\\
\end{cases}
\end{equation}

\begin{equation}\label{eq:ineq2}
\begin{cases}
M_{\gamma_n}(pT,\xi)\ge\gamma_p(pT,\xi) & \forall(n, p)\in N^2\\
M_{\gamma_n}(pT,\chi)\ge\beta_p(pT,\chi) & \forall(n, p)\in N^2\\
M_{\gamma_n}(nT+\frac{\chi-\xi}{v_f},\chi)\ge\beta_p(nT+\frac{\chi-\xi}{v_f},\chi)
& \forall(n, p)\in N^2\\
\hfill \text{s.t.}\quad nT+\frac{\chi-\xi}{v_f}&\in [(p-1)T, pT]
\end{cases}
\end{equation}

\begin{equation}\label{eq:ineq3}
\begin{cases}
M_{\beta_n}(pT,\xi)\ge\gamma_p(pT,\xi) & \forall(n, p)\in N^2\\
M_{\beta_n}(nT+\frac{\xi-\chi}{w},\xi)\ge\gamma_p(nT+\frac{\xi-\chi}{w},\xi) & \forall(n, p)\in N^2\\
\hfill \text{s.t.}\quad nT+\frac{\xi-\chi}{w} &\in [(p-1)T, pT]\\
M_{\beta_n}(pT,\chi)\ge\beta_p(pT,\chi) & \forall(n, p)\in N^2
\end{cases}
\end{equation}
Every single link in a traffic flow control model built upon the Lax-Hopf framework needs to satisfy these constraints. Usually, the control variables are inflows and outflows, in which these constraints are linear. In addition to the compatibility condition for each link, extra constraints are needed for certain control problems. The following sub-section introduces the extra constraints for the VSL and boundary flow control problem. \\
\subsection{Extra Constraints for VSL and Boundary Flow Control}
In a VSL and boundary flow control model, the decision variables are inflows at the mainlane entry node and the on-ramp flows, and the speed limit on certain links. For the links with VSL, the jam density $\rho_m$ and backwave speed $w$ are assumed to be constant while the critical density $\rho_c$ changes with free flow speed $v_f$. Under the triangular FD assumption, we have
\begin{equation}\label{eq:cri_free}
	\rho_c=-\frac{\rho_mw}{v_f-w}
\end{equation}.\\
In addition, the speed limit is assumed to take values only from a discrete set of values that lie within a realistic range $C_{v_f} = \{v_{f_1}, v_{f_2}, ..., v_{f_S}\}$. This assumption will allow us to linearize the non-linear constraints, which is shown later. Let $C_{\rho_c} = \{\rho_{c_1}, \rho_{c_2}, ..., \rho_{c_S}\}$ be the set of corresponding critical densities according to \eqref{eq:cri_free}. Similarly, let $C_{Q} = \{Q_1, Q_2, ..., Q_S\}$ be the set of corresponding capacities. \\
Before introducing extra constraints, the expression of Moskowitz solutions on VSL links need to be modified first. Equations \eqref{eq:mmk}-\eqref{eq:mbeta} are the formulae to compute the Moskowitz solutions under static speed limit. Although the compatibility conditions still hold, Equations \eqref{eq:mmk}-\eqref{eq:mbeta} cannot be readily used under VSL. In reference \cite{vishnoi2021variable}, using Equations \eqref{eq:mmk}-\eqref{eq:mbeta} as the building block, a framework to compute the Moskowitz solutions under VSL was proposed. In short, for the link with VSL, the proposed method divided the temporal space into periods according to the time points where the speed limit changes. Let $T_1, T_2, T_3, ..., T_n$ denote such time points, then the simulation time $[0,t_{\max}]$ is divided into following periods: $[0, T_1], [T_1, T_2], [T_2, T_3], ..., [T_n, t_{\max}]$. Within each period, the Moskowitz function can be computed using \eqref{eq:mmk}-\eqref{eq:mbeta}. Consequently, the densities on the link at the end of a period can be computed as the negative partial derivative of the Moskowitz function with respect to $x$, and these densities can be regarded as the initial densities for the next period. Repeat this process from the first to the last period. For the details, the readers are referred to \cite{vishnoi2021variable}. In the VSL limit control, the modified formulae are substituted into the compatibility conditions \eqref{eq:ineq1}-\eqref{eq:ineq3}. In addition, two different types of constraints and variables are added or modified.
\subsubsection{Non-linear Terms in Compatibility Conditions}
There are two non-linear terms of the decision variables in the compatibility conditions for the links with VSL. The first term is $\rho_cv_f$ in \eqref{eq:mmk} due to the relationship between $\rho_c$ and $v_f$ shown by \eqref{eq:cri_free}. This term can be linearized by
\begin{subequations}\label{eq:linearize1}
	\begin{align}
		\rho_cv_f &= \sum_{s=1}^{S}\delta_s\rho_{c_s}v_{f_s}\\
		\sum_{s=1}^{S}\delta_s &= 1, 
	\end{align}
\end{subequations}\\
where $\delta_s, s\in [1, S]$ are binary variables.\\
The second non-linear term is $q_{in}(i)/v_f$ in \eqref{eq:mgamma}. This term can be linearized by
\begin{subequations}\label{eq:linearize2}
	\begin{align}
	0\le k_{a_s,i}\le \delta_s\frac{Q_s}{v_{f_s}},\forall i,\quad \forall s \label{eq:2a}\\
	\frac{q_{in}(i)}{v_{f_s}}-(1-\delta_s)\frac{Q_{max}}{v_{f_s}}\le k_{a_s,i} \le \frac{q_{in}(i)}{v_{f_s}}, \forall i,\quad \forall s \label{eq:2b}\\
	k_{in,i} = \sum_{s=1}^{S}k_{a_s,i}, \forall i \label{eq:3a}
	\end{align}
\end{subequations}\\
where $\delta_s, s\in [1, S]$ are the same binary variables in \eqref{eq:linearize1}, $Q_{max}$ is the maximum in $C_Q$. $k_{a_s,i}$ are auxiliary continuous variables. \eqref{eq:2a}-\eqref{eq:2b} ensure only one $k_{a_s,i}=\frac{q_{in}(i)}{v_{f_s}}$, and all others equal $0$. Then, we can equivalently replace the non-linear term $q_{in}(i)/v_f$ in the compatibility conditions with $k_{in,i}$.
\subsubsection{Demand and Supply}
Demand and supply are the basic components for the node flow model. Demand is the maximum outflow assuming the downstream capacity is infinity, and supply is the maximum inflow assuming the upstream demand is infinity. Let $L_n$ be the maximal number of vehicles can leave a link up to and including time-step $n$ assuming the downstream capacity is infinity, then the demand for time-step $n$ can be expressed as,
\begin{equation}
	D_n=\frac{L_n-L_0}{T}-\sum_{i=1}^{n-1}q_{out}(i)
\end{equation}
where $L_n$ can be defined through constraints involving the Moskowitz solutions and new auxiliary binary variables. Similarly, the supply can be defined in the same manner. For simplicity, the expression of the constraints are omitted in this paper. Please see \cite{vishnoi2021variable} for details.
\section{Two-stage Modeling Considering Uncertainties in Demand}\label{sec:two-stage}
Using the linear constraints illustrated in Section \ref{sec:review}, Vishnoi and Claudel \cite{vishnoi2021variable} proposed a VSL and boundary flow control model as an MILP in which the demand at the entry nodes are time-dependent but fixed and known. However, in reality, the demand is usually random, and it is of importance to investigate the impact from such uncertainties. Triggered by this motivation, this section tackles this issue using a two-stage stochastic optimization model. The notations used in this paper are summarized in Table \ref{tab:para}. Note that although the inflows and outflows of all links are decision variables in our optimization models, we do not ``control" these flows since they are uniquely determined by the compatibility conditions \eqref{eq:ineq1}-\eqref{eq:ineq3}, the demand and supply constraints shown in previous section after the boundary and VSL control is obtained. Only the speed limit and entry flows are controlled in reality.\\
 \begin{table}[h]
	\caption{Notations}
	\centering
	\begin{tabular}{p{2cm} p{5cm}}
		\Xhline{3\arrayrulewidth}
		\textbf{Notation} & \textbf{Definition} \\
		\Xhline{3\arrayrulewidth}
		\multicolumn{2}{c}{\textbf{Sets}}\\
		\textbf{\textit{Flow Model:}} & \\
		$I$ & set of junctions \\
		$L$ & set of all links\\
		$L_{in}$ & set of entry links\\
		$L_{out}$ & set of exit links\\
		$L_{VSL}$ & set of links with VSL control\\
		$N$ & set of all time steps of the simulation\\
		$N_{T_1}$ & set of all time steps of a project horizon\\
		$N_{T_2}$ & set of all time steps until speed limit is updated\\
		$C_{v_f}^l$ & set of speed limits for link $l$, $\forall l\in L_{VSL}$\\
		$C_{\rho_c}^l$ & set of critical densities for link $l$, $\forall l\in L_{VSL}$\\
		$C_Q^l$ & set of capacities for link $l$, $\forall l\in L_{VSL}$\\
		$C_D$ & set of possible joint demand for all entry links\\
		\Xhline{4\arrayrulewidth}
		\multicolumn{2}{c}{\textbf{Parameters}}\\
		\textbf{\textit{Flow Model:}} & \\
		$\psi$ & fundamental diagram\\
		$v_f^l$ & free flow speed of link $l$ (m/s), $\forall l\in L\setminus L_{in}$\\
		$\rho_c^l$ & critical density of link $l$ (veh/m), $\forall l\in L\setminus L_{in}$\\
		$w^l$ & backwave speed of link $l$ (m/s), $\forall l\in L\setminus L_{in}$\\
		$\rho_m^l$ & jam density of link $l$ (veh/m), $\forall l\in L\setminus L_{in}$\\
		$Q^l$ & capacity of link $l$ (veh/s), $\forall l\in L\setminus L_{in}$\\
		$\tilde{d(l)}$ & random demand for entry link $l$ (veh/s),  $\forall l\in L_{in}$\\
		$P$ & a vector of probability of each scenario in $C_D$\\
		\midrule
		\textbf{\textit{Discretization:}} & \\
		$X$ & length of spatial segment (m)\\
		$T$ & time step size (s)\\
		$k_{\text{max}}(l)$ & number of spatial segments on link $l$\\
		$n(.)$ & number of time steps in a horizon\\
		\Xhline{4\arrayrulewidth}
		\multicolumn{2}{c}{\textbf{Control variables}}\\
		$q_{\text{in}}(l,t)$ & inflow of link $l$ at time $t$ (veh/s)\\
		$q_{\text{out}}(l,t)$ & outflow of link $l$ at time $t$ (veh/s)\\
		$q'_{\text{in}}(l,t)$ & boundary control for entry link $l$ at time $t$ (veh/s)\\
		$\delta_j^l$ & binary variable: 1, if the $j$th element in $C_{v_f}^l$ is the optimal speed; 0, otherwise\\
		\Xhline{3\arrayrulewidth}
	\end{tabular}\label{tab:para}
\end{table}\\
\subsection{Two-stage Optimization Model}\label{sec:two_stage}
\begin{figure}[htb!]
	\includegraphics[width=3.5in]{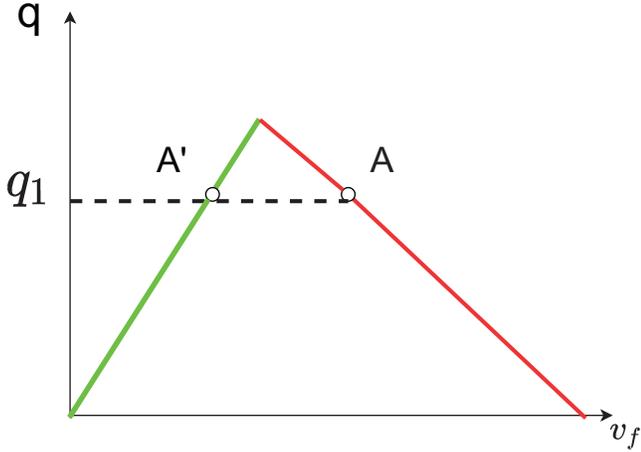}
	\centering
	\caption{Mitigate congestion while maintaining mobility.}
	\label{fig:tfd}
\end{figure}

The basic idea for a two-stage model is that the optimal decisions at the first stage, which is called a ``here-and-now" decision should be made before the realization of random parameters being observed. After the uncertain data is revealed, the optimal solution at the second stage, which is called a ``recourse" decision, is obtained by solving a recourse optimization model in which the realization of the random variable and the optimal solution from the first stage are model parameters. The two-stage optimization model ensures the optimal control from the first stage leads to a combined optimal value of the object function at the first stage and the expectation of the recourse function at the second stage.\\
In the context of control problem in this paper, the demands of the entry links are random, we have to make decisions on how many vehicles are allowed to entry the road section (the first stage) before we observe the real demand. After implementing the boundary control and observing the real demand, the speed limits on certain downstream links are optimized (the second stage) for certain objectives. We assume the joint distribution of the demands $P$ is discrete and known. Then, the two-stage model can be expressed as
\begin{equation}\label{eq:twostage}
\begin{split}
\max_{\bm{q'_{in}}} \quad & \omega_0\sum_{l\in L_{in}}\sum_{t=1}^{N_{T_1}}q'_{in}(l,t)+\mathbb{E}[g(\bm{q'_{in}},\bm{\tilde{d}})]\\
\text{s.t.} \quad & \bm{q'_{in}} \in \bm{R} \\
& \text{OCC}, \quad \forall l\in L\\
& \text{DSC}, \quad \forall i\in N\\
\end{split}
\end{equation}
The first stage aims to find the optimal boundary flow control $q'_{in}(t,l), \quad \forall t\in N, l\in L_{VSL}$, i.e., the maximal inflows that are allowed at each time step for the incoming boundary links, considering the VSL control effect in the second stage. OCC and DSC stand for the original compatibility condition and demand-supply constraints at junctions, respectively. The goal of the first stage is to find the feasible set for the boundary flow control that satisfy these constraints without optimizing specific objective functions. The objective function of the first stage is the maximization of weighted sum of boundary controls in which $\omega_0$ is a small weight. The reason of this term will be explained after we finish the second stage. $\mathbb{E}[g(\bm{q'_{in}},\bm{\tilde{d}})]$ is the expectation of the objective value of the second stage model, in which the optimal solution from the first stage and the realization of demand are model parameters.\\
After we implement the first stage control $\bm{q'_{in}}$, a realized scenario of the joint demands, $\bm{d_j}$, can be observed. Then, the recourse function in the second stage can be expressed as
\begin{align}\label{eq:2nd}
\begin{split}
& g(\bm{q'_{in}},\bm{d_j}) =\\
& \max_{\bm{v_j^f}} \quad \!\begin{aligned}[t]
& \sum_{l\in L_{out}}\sum_{t=1}^{N_{T_1}}q_j^{out}(l,t)(N_{T_1}-t+1) \\
& - \sum_{l\in L_{VSL}}\sum_{t=1}^{N_{T_1}}\omega_1^lq_j^{in}(l,t) \\
& - \sum_{l\in L_{in}}\sum_{t=1}^{N_{T_1}}\omega_2^lq_j^{in}(l,t) \\
& - \sum_{l\in L_{in}}\sum_{t=1}^{N_{T_1}}\omega_3^l\left(\sum_{t_1=1}^td_j(l,t_1) - \sum_{t_1=1}^tq_j^{in}(l,t_1)\right)(1+e)\\
& - \sum_{l\in L_{in}}\sum_{t=1}^{N_{T_1-1}}\omega_4^l|q_j^{in}(l,t) - q_j^{in}(l,t+1)|
\end{aligned}\\
& \text{s.t.} \quad q_j^{in}(l,t)\le q'_{in}(l,t) \quad \forall l \in L_{in}, t\in N_{T_1}\\
& \sum_{t=1}^{t_1}q_j^{in}(l,t)\le \sum_{t=1}^{t_1}d_j(l,t), \quad \forall l\in L_{in},t_1\in  N_{T_1}\\
& q_j^{in}(l,t) + M\delta_j(l,t) \ge q'_{in}(l,t), \quad \forall i\in L_{in}, t\in  N_{T_1}\\
& \sum_{t=1}^{t_1}q_j^{in}(l,t)+M[1-\delta_j(l,t)]\ge \sum_{t=1}^{t_1}d_j(l,t), \quad \forall l\in L_{in},t_1\in  N_{T_1}\\
& \delta_j(l,t) \in \{0,1\}, \quad \forall l\in L_{in}, t\in  N_{T_1}\\
& \text{OCC}, \quad \forall l\in L\setminus L_{VSL}\\
& \text{MCC}, \quad \forall l\in L_{VSL}\\
& \text{DSC}, \quad \forall n\in N\\
\end{split}
\end{align}
where $j$ denotes the index of scenarios of the joint demands, $q^{in}_j(l,t)$ and $q^{out}_j(l,t)$ indicates the inflow and  outflow of link $l$ at time $t$. The control variables in the second stage are the speed limits for all links with VSL control, denoted by $v^f_j(l), \forall l\in L_{VSL}$. The first term in the objective function is to maximize the sum of weighted of outflows over all exit links. The weight, which is a decreasing function of time step, is used to force vehicles to proceed as long as there is available space ahead so that we can avoid unnecessary stops. The second term is to minimize the weighted inflows of VSL links. This term is to limit inflows on VSL links when congestion takes place downstream. $\omega_1^l$ is a small weight to ensure the flow reduction on these VSL links will not sacrifice the overall throughput of the road section, which is the first term in the objective function. For the same reason, the third term is used to minimize the weighted inflows of the entry links. Through using this term, the traffic condition, in terms of level of congestion, on the whole road section is controlled. Again, the weight should be small enough to not sacrifice the throughput. Figure \ref{fig:tfd} shows a triangular FD for a single link to illustrate this idea. Assume at the beginning of the simulation, both the supply of downstream node and the demand of upstream node are equal to the capacity. Regardless of the initial density, after a certain period, this link will reach to the stable state where both its inflow and outflow are equal to capacity. Let us assume the downstream capacity drops to $q_1$ at some point during the simulation due to a fixed bottleneck, without these two terms, to maximize the outflows, the control may still send vehicles onto this link at its capacity until the backward-moving shockwave arrives the entry node. Then, both the inflows and outflows afterwards equal $q_1$, and the traffic on this link will be stablized at point A which belongs to the congestion phase of the FD. However, by introducing the second term, when the capacity drop happens downstream, the control will reduce inflows to some extend at the beginning as long as it does not sacrifice outflows. As a result, the density on the link will decrease until it reaches the density corresponding to point A' which belongs to the free flow phase. Therefore, this term can mitigate congestion and maintain the outflows concurrently. The fourth term is a penalty term for blocking vehicles at the entry nodes. $\left(\sum_{t_1=1}^td_j(l,t_1) - \sum_{t_1=1}^tq_j^{in}(l,t_1)\right)$ is the cumulative sum of blocked vehicles at time $t$, and $e$ is the initial queue length at the beginning of the project horizon, which implies a higher penalty from a longer initial queue. Therefore, there is a trade-off between the third and the fourth term. Lastly, the objective for the fifth term, which is the minimization of the difference between consecutive inflows at the entry nodes, is twofold: first, since safety is affected by the temporal variation in traffic flow \cite{golob2004freeway, theofilatos2014review, yuan2021developing}, this term aims to smooth the entry flows; second, without this term, the control might send very few or even stop sending vehicles in the last few steps as long as there are enough vehicles leaving the exit nodes and the third term surpasses the fourth penalty term, which could block a large number of vehicles and cause a very large penalty weight for the next horizon.\\
The first constraint means the inflows cannot exceed the control from the first stage. The second constraint indicates at any time step $t$, the cumulative number of inflows is upper bounded by the cumulative demand. The third and fourth constraints force the inflow at each time step is equal to the maximal value satisfying the first two constraints, where $\delta_j(i,t)$'s are binary variables and $M$ is a big value. MCC is the modified compatibility conditions for links with VSL control that are explained in the previous section.\\
Now, we explain the objective function of the first stage in \eqref{eq:twostage}. In short, the goal of this term is to ensure uniqueness and practicability. Let us assume the optimal control from a deterministic model is equal to the assumed demand $d_1$. Since this control allows all vehicles to enter the road, any control with a higher value will lead to the same traffic operation and objective value. To make the control unique, by adding the term in the first stage with a small weight $\omega_0$, we force the boundary control to be equal to the capacity without changing the control performance. Note that this term only affects the cases in which all demands are allowed by the optimal control. If the optimal control is less than demand, increasing optimal control will change the inflows and lose the optimality. In addition to the mathematical consideration, we think this approach also makes practical sense. If a control is equal to the demand, all vehicles are allowed to enter, and thus, it is equivalent to "no control". Additionally, we will implement the control from both the proposed two-stage model and the deterministic model on scenarios with random demands and compare their objective values. Using this term can also prevent the deterministic control from blocking vehicles if the real demand is higher than the assumed demand.\\
For random demands with discrete distributions, the expectation of the recourse function can be expressed as
\begin{equation}\label{eq:milp}
\mathbb{E}[g(\bm{q'_{in}},\bm{\tilde{d}})] =\sum_{j=1}^{|C_D|} p_jg(\bm{q'_{in}},\bm{d_j})
\end{equation} 
Above all, Equation \eqref{eq:twostage} is an MILP with the boundary inflows and speed limit for each scenario as decision variables. The number of decision variables and constraints are linear to $J$. 
\subsection{Rolling Horizon Scheme}
\begin{figure}[htb!]
	\includegraphics[width=3.5in]{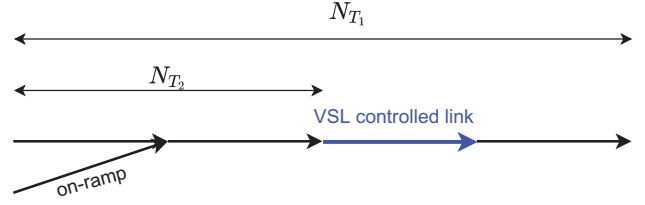}
	\centering
	\caption{An example corridor.}
	\label{fig:corridor}
\end{figure}

\begin{figure}[htb!]
	\includegraphics[width=3.5in]{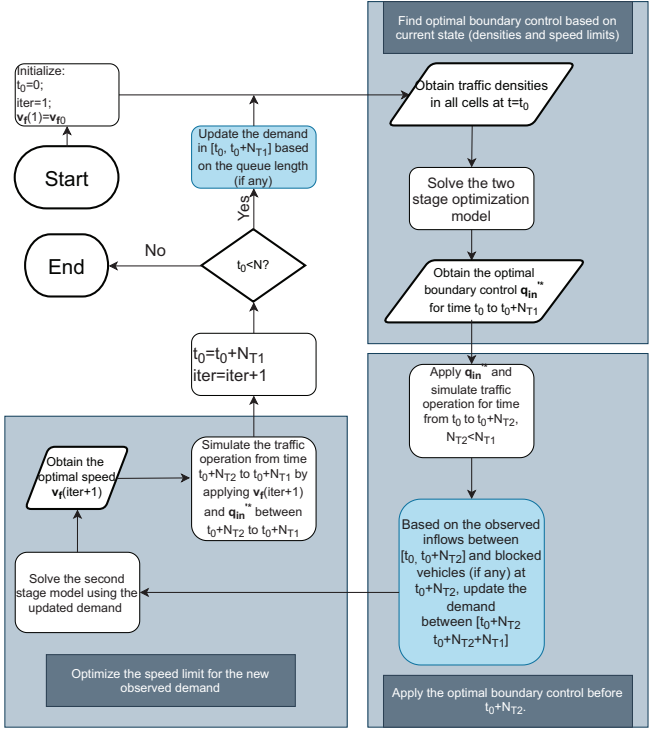}
	\centering
	\caption{Rolling horizon scheme.}
	\label{fig:rhs}
\end{figure}

We utilize a rolling horizon scheme to coordinate the boundary and VSL control. Let us consider a highway corridor shown in Figure \ref{fig:corridor}. Let $N_t$ denote the time needed to finish a trip on this corridor with free flow speed, i.e., $N_{T_1} = R/(v_fT)$ where $R$ is the length of the corridor, and $N_{T_2}$ denote the free flow travel time from the entry node of the corridor to the upstream node of the first link with VSL control. At the beginning of a project horizon $t=t_0$, the optimal boundary control is solved for the time duration $[t_0, t_0+N_t]$. Note we assume the demand does not change with time in $N_{T_1}$ although the demand is random. Then, only the control before $t_0+N_{T_2}$ is implemented. Then, we observe the real demand in this project horizon, and VSL are updated at $t_0+N_{T_2}$. $N_{T_1}$ is called project horizon, and $N_{T_2}$ is referred to as rolling horizon. As a result, the speed limit is updated after the real demand is observed and before the first vehicle arrive at the first VSL link.\\
\\
Figure \ref{fig:rhs} provides a detailed explanation on the iterative procedures. Note that there are two places, shown in blue in Figure \ref{fig:rhs}, where the demand need to be updated based on the observed queue. Let $e$ be the queue length in terms of flow. The update at the beginning of an iteration (top middle in Figure \ref{fig:rhs}) includes following steps:\\
\textit{Step 1:} Initialize random demand matrix $d_{ini}$ with size $N_{T_1} \times |C_D|$. Each column is a constant vector with values equal to the element in $C_D$.\\
\textit{Step 2:} Update each column of $d_{ini}$. We use the first column to explain the process. Let $j=1$ and $e'=\min\{Q-d_{ini}(j,1), e\}$. Then, update $d_{ini}(j,1) = d_{ini}(j,1)+e'$, $e=e-e'$ and $j=j+1$. Repeat this process until $j>N_{T_1}$ or $e=0$. Repeat this process for all columns.\\
\textit{Step 3:} Use the updated demand for the optimization model, i.e., $\bm{d_j}$ in Equation \eqref{eq:milp} is the $j$th column in the updated demand matrix.\\
The demand can be updated in a similar way at $t_0+N_{T_2}$ except for the initialization. At $t_0+N_{T_2}$, we observe the demand in the current project horizon, so the initialized vector ($N_{T_1}\times 1$) has the first $N_{T_1} - N_{T_2}$ elements equal to the observed demand, and the rest elements equal to the mean value of the demand. Remember we assume the demand during one project horizon does not change, and the traffic operation is simulated based on the Lax-Hopf solution, which has been utilized in \cite{vishnoi2021variable, liu2021robust}.

\section{Case Study}\label{sec:casestudy}
\begin{figure}[htb!]
	\includegraphics[width=3in]{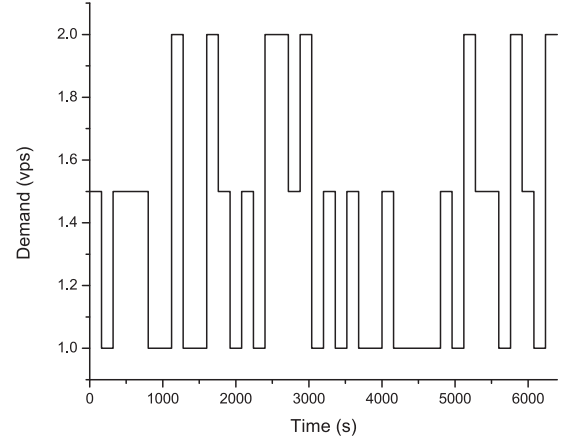}
	\centering
	\caption{An example for demand evolution.}
	\label{fig:demand}
\end{figure}

In this section, we apply the proposed model and the deterministic models with different demand input, i.e., mean value, minimum and maximum to the highway corridor shown in Figure \ref{fig:corridor} and compare their objective values.\\
This corridor consists of 4 main links and one on-ramp, and each link has 4 lanes. The third link is VSL controlled, and other links have fixed speed limit. All links are 1200 m and are divided into 2 even segments. Also, a total simulation time of 6400 s is divided into 320 even segments, so each time step is 20 s. For any link $l \in L\setminus L_{in}$, the free flow speed $v_f = 30$ m/s (67 mph), the critical density $\rho_c=0.0175$ veh/m/lane (28 veh/mile/lane), the jam density $\rho_m = 0.125$ veh/m/lane (201 veh/mile/lane) and the backwave speed $w=-4.9$ m/s (-11 mph). Therefore, the capacity is equal to 2.1 vps (1890 vphpl). The set of speed limits is $C_{v_f} = \{10, 15, 20, 25, 30\}$. Each project horizon $N_{T_1}$ has 8 time steps (160 s), and the rolling horizon $N_{T_2}$ has 4 time steps (80 s). Let the demand at the entry node be random and have a discrete distribution. Let $C_D = \{1, 1.5, 2\}$ (vps), and the corresponding distribution $P = \{0.4, 0.2, 0.4\}$. Assume at time $t=0$, the downstream capacity drops to $1.4$ vps due to certain reasons such as a crash, so we implement the model \eqref{eq:twostage} to control the boundary flow in order to reduce congestion and maintain safety. In addition, we assume that the drivers on the frontage road are aware of this situation so the travel demand on the on-ramp is low and equals 0.05 vps. Consequently, we assume at the merging node, the vehicles from on-ramp can always enter as long as there is available space on the downstream link. This can be modeled with linear constraints in a similar way as the demand-supply constraints. Note this assumption is for simplicity, and other node models can also be used. The weights in the objection function are $\omega_0=1e^{-4}$, $\omega_1=0.01$, $\omega_2=0.02$, $\omega_3=0.003$ and $\omega_4=10$. A large value of $\omega_4$ is used to ensure smooth entry flows.\\
Following the demand distribution, we generate ten scenarios with different random seeds at the entry node. For each scenario, we apply the proposed two-stage model, a deterministic model with minimal demand, a deterministic model with mean demand, and a deterministic model with maximal demand and compare their objective values. From here on, these four models are referred to as D-Min, D-Mean, D-Max and Two-stage, respectively. To provide an example, Figure \ref{fig:demand} shows the demand of the second scenario.

\subsection{Results}
\subsubsection{Non-comparable terms}
\begin{figure}[htb!]
	\includegraphics[width=3in]{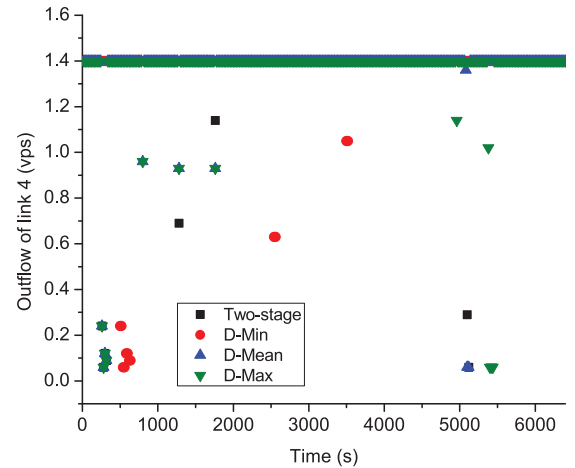}
	\centering
	\caption{Comparison between outflows of the corridor.}
	\label{fig:outflow}
\end{figure}

To validate the proposed model, the control performance in terms of the objective value during the whole simulation are compared. However, there is no significant difference in the total objective value among these three model, which is expected. The reason is that the maximization of downstream outflow has the highest priority in our objective function, the first term in Equation \eqref{eq:2nd} has much larger weights than other terms. In addition, the mean value of demand (1.5 vps) is higher than the downstream capacity (1.4 vps), so the downstream outflows in all models are very similar and close to the capacity. As a result, the objective value of this term is very close, and the differences in other terms are obscured. Figure \ref{fig:outflow} shows the comparison for the 2nd scenario.\\
Two-stage has higher downstream flow in 1 scenario and lower downstream flow in 3 scenarios than D-Min, and both models have the same value in 6 scenarios. It has higher downstream flow in 3 scenarios and lower downstream flow in 1 scenario than D-Mean, and both models have the same value in 6 scenarios. It has higher downstream flow in 9 scenarios and lower downstream flow in 1 scenario than D-Max. However, as mentioned above and shown in Figure \ref{fig:outflow}, the difference is negligible. In addition to this term, the differences in terms associated with $\omega_1$ and $\omega_2$ are also negligible. In all these four models, the VSL control is performed after observing the real demand, so intuitively, all these controls are capable of sending enough vehicles to link 3 to maximize the downstream outflow without causing congestion on links 3 and 4, so they generate close values for the term associated with $\omega_1$. Therefore, although the inflows to link 3 in single project horizon are different among the three controls, the sum over all project horizons are almost equal. Similarly, the difference in the sum of term associated with $\omega_2$ is indistinct as well.\\
\subsubsection{Comparable terms}
\begin{figure}[!t]
	\centering
	\subfloat[Block penalty]{\includegraphics[width=2.5in]{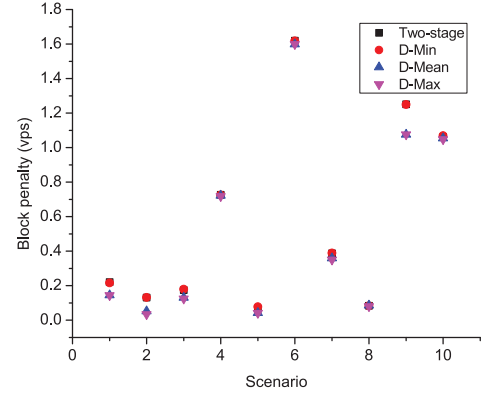}
		\label{Block penalty}}
	\hfil
	\subfloat[Inflow fluctuation]{\includegraphics[width=2.5in]{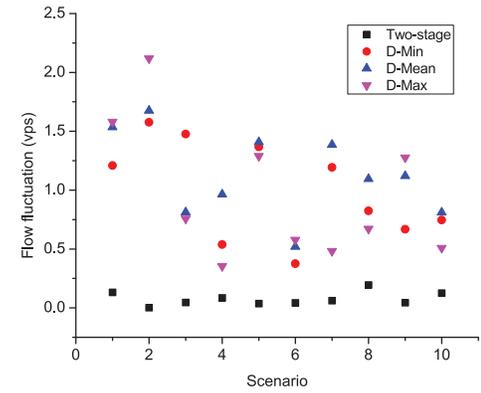}
		\label{Inflow fluctuation}}
	\hfil
	\subfloat[Combined value]{\includegraphics[width=2.5in]{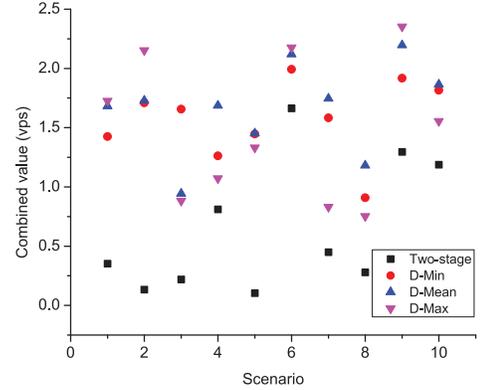}
		\label{sum}}
	\caption{Comparison of objective values.}
	\label{fig:somparison_ov}
\end{figure}

\begin{figure}[!t]
	\centering
	\subfloat[Two-stage]{\includegraphics[width=3in]{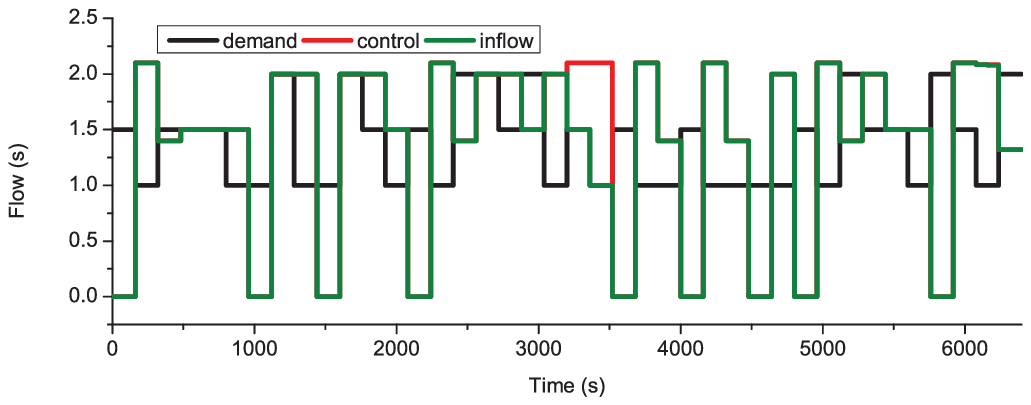}
		\label{fig:two-stage}}
	\hfil
	\subfloat[D-Min]{\includegraphics[width=3in]{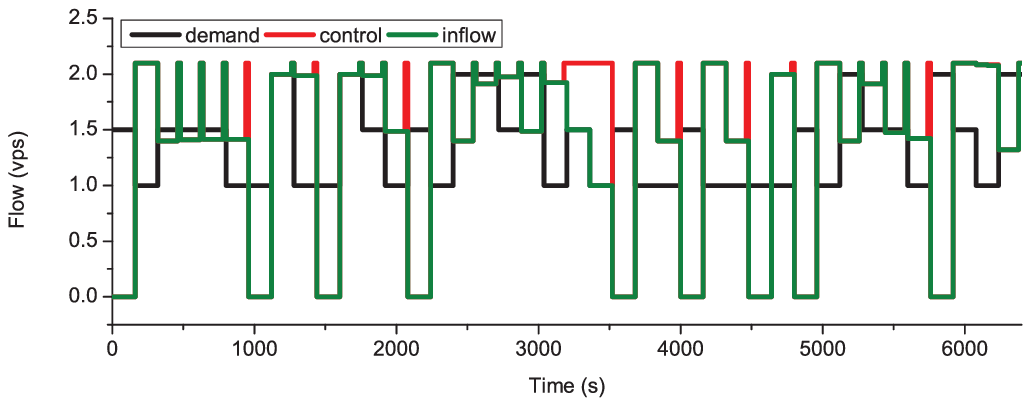}
		\label{fig:D-min}}
	\hfil
	\subfloat[D-Mean]{\includegraphics[width=3in]{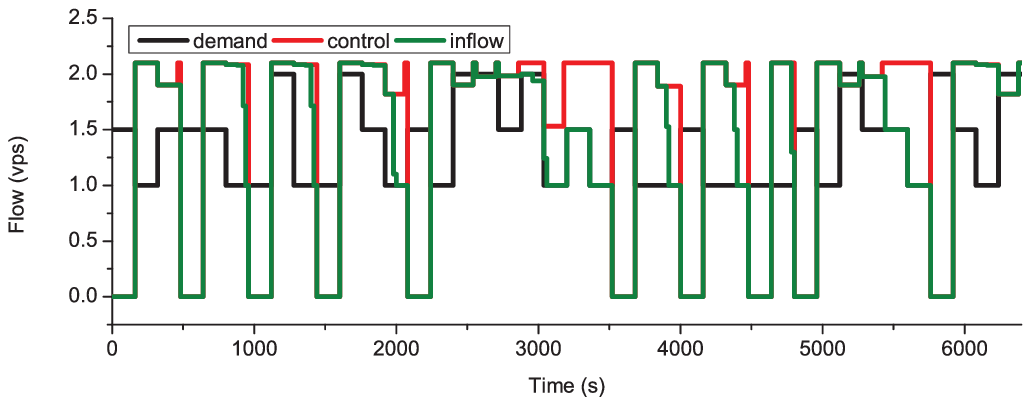}
		\label{fig:D-Mean}}
	\hfil
	\subfloat[D-Max]{\includegraphics[width=3in]{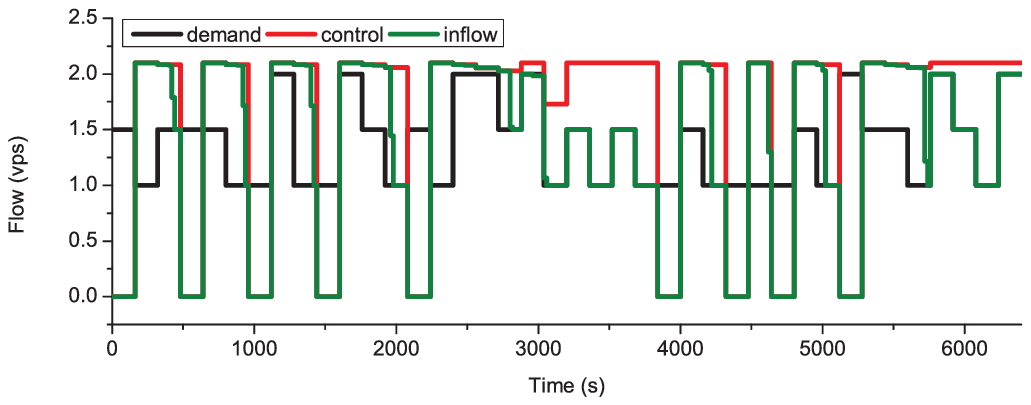}
		\label{fig:D-Max}}
	\caption{Comparison of traffic evolution.}
	\label{fig:somparison_evo}
\end{figure}

\begin{figure}[htb!]
	\includegraphics[width=3.5in]{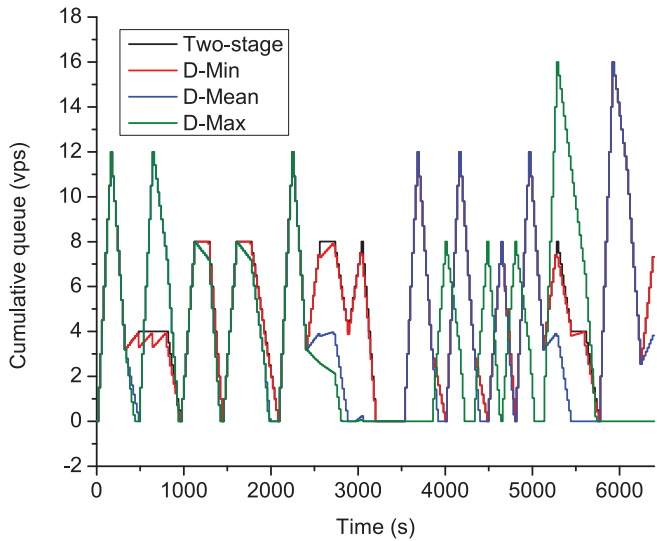}
	\centering
	\caption{Comparison of cumulative queue length.}
	\label{fig:queue}
\end{figure}

\begin{figure}[htb!]
	\includegraphics[width=3.5in]{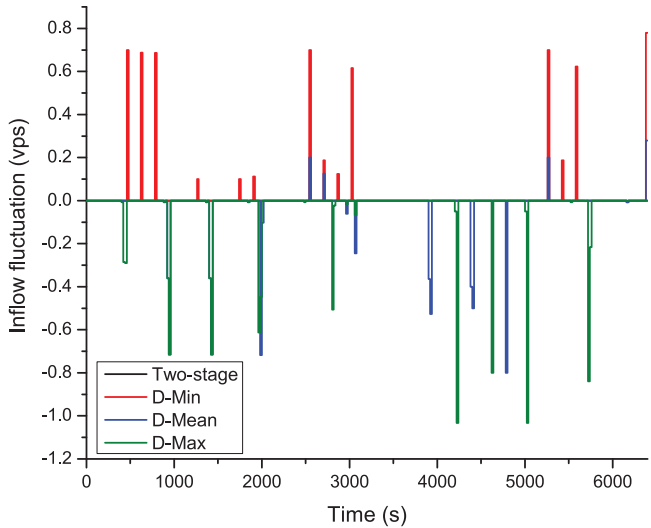}
	\centering
	\caption{Comparison of temporal inflow fluctuation.}
	\label{fig:compare_fluc}
\end{figure}

Unlike the terms mentioned above, the value of last two terms, which are called block penalty and flow fluctuation, in different project horizons cannot compensate each other. We focus on these two terms in the rest of this paper. Figure \ref{fig:somparison_ov} shows the comparison of these two terms and their sum. It is seen there is a trade-off between the block penalty and flow fluctuation. The proposed model and D-Min have higher block penalty and lower flow fluctuation than the other two models. In addition, the flow fluctuation from Two-stage is considerably lower than the other three models. Overall, the proposed model outperforms all deterministic models in the combined value. The proposed model can reduce the sum of both terms by 58.7\%, 60.9\% and 56.2\% compared to D-Min, D-Mean and D-Max, respectively. \\
\\
To provide more detail, Figure \ref{fig:somparison_evo} shows the evolution of demand, same as Figure \ref{fig:demand}, boundary control and real inflow of a scenario. The pattern of control and inflow from these models are significantly different. Note demand and inflow are instantaneous rather than cumulative, so the inflow at certain time points could be higher than the demand, as shown in the plots, if there are vehicles blocked from previous time steps. Although the control within most project horizons from all models is constant, Figure \ref{fig:somparison_evo} clearly shows that the discrepancies between boundary control and real inflow in deterministic models occurs much more frequently than Two-stage. This discrepancy could result from both overestimation and underestimation of the demand, and it leads to significant flow fluctuations. It shows that Two-stage is able to reduce the occurrence of such discrepancies. As a result, it has a much lower flow fluctuation penalty, as shown in Figure \ref{Inflow fluctuation}. The reason for fluctuation resulting from deterministic models is explained in greater detail later in this section. 
\\
Figure \ref{fig:queue} shows the comparison of cumulative queue lengths. Note the cumulative queue has a unit of vps, and it can be multiplied by the time step size to obtain the cumulative number of blocked vehicles. Since the blocked penalty term in \eqref{eq:2nd} is positively related to the queue length at the beginning of a horizon, the optimization models will dissipate the blocked queue if it exceeds a certain value. This is the reason of peaks on the curves. Also, it shows overall, the queues resulting from D-Mean and D-Max are shorter than Two-stage and D-Min, which is consistent with the 2nd scenario in Figure \ref{Block penalty}. This is because D-Min underestimates the real demand and generates control lower than the real demand, while Two-stage tries to reduce inflow fluctuation through limiting the number of entering vehicles. \\
Figure \ref{fig:compare_fluc} shows a more straightforward comparison of inflow fluctuation, in which each point is the difference of inflows between consecutive time steps. Note that only the fluctuation within each project horizon is penalized, and we do not consider fluctuations at connections of project horizons. Since the demand is assumed to follow a discrete distribution, and in each project horizon, the demand is a constant drawn from the distribution, the demand fluctuation at the connection of consecutive project horizons is inherent. The boundary control is adjusted according to the change of demand, so the inflow fluctuation at the connection points are inevitable, as shown in Figure \ref{fig:somparison_evo}. Therefore, Figure \ref{fig:compare_fluc} only shows the fluctuation within each project horizon, which is also the parameter the objective function aims to minimize.\\ While the flow variations in all deterministic models are significant, the proposed model generates highly smooth flow with only one non-zero (0.0078 vps) inflow fluctuation in this example. The flow fluctuation can be divided into two types corresponding to the generation mechanism. The first type of fluctuation results from the overestimation of demand, and the value of such fluctuations is negative. This type exists in both D-Mean and D-Max. The reason is as follows. Overall, D-Max overestimates the real demand, so its boundary control can also be higher than the real demand. Let us consider project horizons with an initial queue and low demand $D_1$, and control $q>D_1$ corresponding to the maximal demand $D_{\text{max}}$. Under this situation, the inflows at the beginning steps are equal to the control $q$ due to the existence of a queue. However, since the real demand is lower than assumed value, the queue can be cleared in the middle of this project horizon. Then, the inflows drop to the demand $D_1$ from $q$, which leads to negative flow fluctuation. Since D-Max always overestimates the demand (except for the accurate estimation), its fluctuation values are all negative. The second type of fluctuation originates from the underestimation of demand, and the values are positive. This type exists in both D-Min and D-Mean. Contrary to D-Max which always overestimates the demand, D-Min can only underestimate the demand. Let us assume the real demand $D_1$ is higher than the assumed value $D_{\text{min}}$, and the control during the last steps equal to $D_{\text{min}}$. As explained in Section \ref{sec:two_stage}, under this situation, to ensure the uniqueness of the optimal control, we force the boundary control to be equal to the capacity. Consequently, the inflows in these steps can possibly increase by certain values depending on the real demand $D_1$. For example, if $D_1$ equals capacity, the inflows will increase to capacity. This is the reason of the positive fluctuation values in D-Min. Since D-Mean can both overestimate and underestimate the demand, both types of fluctuation exist in D-Mean, as shown in Figure \ref{fig:compare_fluc}. The fluctuations in Figure \ref{fig:compare_fluc} are also consistent with the evolution of inflows in Figure \ref{fig:somparison_evo}. Since the proposed model takes all demand scenarios and the corresponding probability distribution into consideration, its boundary control outperforms the deterministic models in a random environment.  
\subsection{Control performance under different demand variation}
\begin{figure}[htb!]
	\includegraphics[width=3.5in]{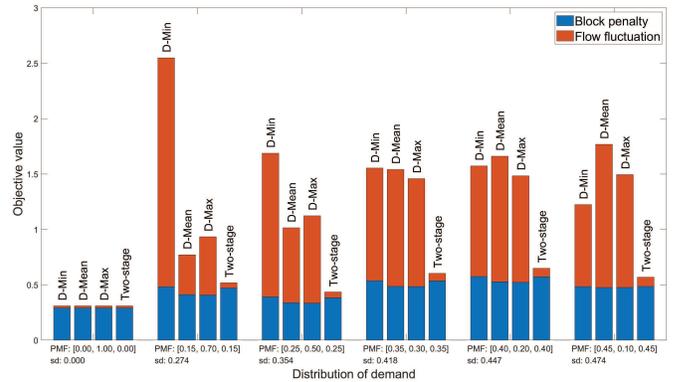}
	\centering
	\caption{Performance comparison with different demand distribution.}
	\label{fig:compare_sd}
\end{figure}

To validate the performance of the proposed model under various demand distributions, this section compares the block penalty term and flow fluctuation term with different demand distributions, as shown in Figure \ref{fig:compare_sd}. In this section, the set of demand values is the same as the previous section, and only symmetric distributions are considered, i.e., the probabilities of $\tilde{d} = 1$ and  $\tilde{d} = 2$ are the same. The label on the horizontal axis represents both the probability mass function (PMF) and the corresponding standard deviation (sd). As expected, when $\sigma_d=0$, all four models have the same model inputs and control performance. When the standard deviation is less than 0.4, D-Mean performs better than D-Min and D-Max since the probability of the mean value, which is used as the input in D-Mean, is relatively high when the standard deviation is small and thus the probability of accurate estimation is high. With the increase of the variation, the probability of the mean value becomes smaller, and D-Min ($sd\ge0.447$) and D-Max ($sd\ge0.418$) outperform D-Mean when $sd>0.4$. For all scenarios with $sd>0$, the proposed model conquers the deterministic models with considerable reduction in the value of flow fluctuation and a slightly higher value of the block penalty term. In addition, D-min suffers from much higher fluctuation than D-Max when sd is small. As mentioned earlier, D-Min increases flow from a smallest value to capacity which could be larger than the flow drop from maximal demand to a slightly lower demand. This gap is mitigated with the increase of the standard deviation. Another observation is that the change of the flow fluctuation value from both deterministic models is much more significant than the block penalty, which is expected. The reason is that the term with the highest weight in the objective function is to maximize the outflow of the downstream link. Consequently, enough vehicles will be allowed to enter the corridor to ensure the downstream link receives enough demand, and the inflows should be similar regardless of the variation in demand.  Also, the mean value of the demand is unchanged among these scenarios, so the number of blocked vehicles should be similar as well. As a result, the change of block penalty is insignificant. Unlike the block penalty, the uncertainties in demand has a considerable impact on flow fluctuation, and the proposed model is capable of smoothing the entry flows while blocking slightly more vehicles and generates much smaller combined objective values. 


\section{Summary and Future Work}\label{sec:summary}
This paper proposes a two-stage stochastic optimization model derived from the Lax-Hopf solution to the LWR model to consider the effect of uncertain demand on boundary and VSL control. In the first stage, based on the current state, the boundary control is calculated. Then, after the real demand is observed, the VSL is obtained in the second stage. The proposed model coordinates the boundary control and VSL control by a rolling horizon scheme. The proposed model is validated under the situation when the downstream capacity is dropped and the demand is higher than the capacity. The control performance is compared to deterministic models with demand inputs equal to the minimum, mean value and maximum of the real demand. As expected, all these models can ensure the maximization of throughput since the demand is higher than capacity and it has a large weight in the objective function. However, the proposed model can reduce the boundary inflow significantly, which plays an important role on safety. Although more vehicles are blocked in the proposed model, the overall objective value is improved significantly. The efficacy of the proposed model is demonstrated with different variations in the demand.\\
The proposed model can be improved in the following ways. First, the Lax-Hopf used in this model assumes vehicles can accelerate instantly, which has a significant impact on both boundary and VSL control. Qiu et al. \cite{qiu2013exact} proposed a Lax-Hopf solution with bounded acceleration. Implementing the bounded acceleration on our control algorithm is a promising research topic to improve the practicability of the proposed model. Second, the computational burden needs to be resolved. The case studies in this paper consider 5 links including on-ramp, 5 scenarios for speed limit and 3 scenarios for demand. There are 8 time steps (160 s) in each project horizon. On average, while it takes less than 1 s to solve the deterministic models, it takes about 15 s to solve the MILP \eqref{eq:twostage} using CPLEX V12.9 in Matlab R2018b on a Microsoft Surface with i5 CPU and 8.0 GB RAM. How to improve the computation speed is another research direction. Third, the demand uncertainty in each project horizon is assumed to be independent and identically distributed (i.i.d), which ignores the correlation between consecutive project horizons and the demand can have unreasonable fluctuations. One possible way to address this issue is to assume the distribution of the current project horizon only depends on the traffic state and control in previous horizon, and model the problem as a Markov Decision Processes (MDP) and find the corresponding control policy.    
\bibliographystyle{IEEEtran}
\bibliography{IEEEabrv,twostage}
\vspace*{-2\baselineskip}
\end{document}